\documentclass[12pt,a4paper]{article}

\usepackage{graphicx}          
\usepackage[dvips]{epsfig}    
\usepackage{url}
\usepackage{amsmath}
\usepackage{amssymb}
\usepackage{amsthm}

\newtheorem{remark}{Remark}
\newtheorem{defi}{Definition}
\newtheorem{example}{Example}
\newtheorem{thm}{Theorem}

\newtheorem{lem}{Lemma}

\newcommand{\re}{{\mathbb R}}

\newcommand{\n}{{\mathbb N}}

\newcommand{\vectun}{{e }}
\providecommand{\com[2]}{\begin{tt}[#1: #2]\end{tt}}

\newcommand{\fulllanguage}{{\tt full language problem}}
\newcommand{\setmat}{{\Sigma{}}}
\newcommand{\conicnorm}{Entrywise-comparison{}}
\title{When is a set of LMIs a sufficient condition for stability?} 

\author{Amir Ali Ahmadi, Rapha\"{e}l M. Jungers, Pablo A. Parrilo,\\ and Mardavij
Roozbehani\thanks{A. A., P. P., and M. R. are with the Laboratory for
Information and Decision Systems, Massachusetts Institute of
Technology.  R. J. is with the F.R.S.-FNRS and with the ICTEAM Institute, Universit\'e catholique de Louvain, Belgium. Email: \{\texttt{a\_a\_a}, \texttt{parrilo},
\texttt{mardavij}\}\texttt{@mit.edu},
\texttt{raphael.jungers@uclouvain.be.}} }        

\begin{document}
\maketitle                       
\begin{abstract}      We study stability criteria for discrete time switching systems.  We investigate the structure of sets of LMIs that are a sufficient condition for stability (i.e., such that any switching system which satisfies these LMIs is stable).  We provide an exact characterization of these sets.  As a corollary, we show that it is PSPACE-complete to recognize whether a particular set of LMIs implies the stability of a switching system.
\end{abstract}

\section{Introduction}  In many practical engineering situations, the dynamical behaviour of the system at stake can be modeled as a switching system like the one represented in Equation (\ref{eq-switching}):
\begin{equation}\label{eq-switching}
x_{k+1}=A_{\sigma\left(k\right)  }x_{k},
\end{equation} 
where
$\setmat\mathrel{\mathop:}=\left\{ A_{1},...,A_{m}\right\} $ is a set of matrices, and the function $$\sigma(\cdot): \n\rightarrow \{1,\dots, m\} $$ is called the \emph{switching signal}.
As a few examples, applications ranging from Viral Disease Treatment optimization (\cite{hmcb10}) to Multi-hop networks control (\cite{pappas-multihop}), or trackability of autonomous agents in sensor networks (\cite{cresp}) have been modeled with switching systems. 

 One of the central problems in the study of switching systems is their stability: do all the trajectories $x(t)$ tend to zero when $t\rightarrow \infty,$ whatever switching law $\sigma\left(k\right)$ occurs?  The answer is given by the \emph{joint spectral radius} of the set $\setmat$ which is defined
as

\begin{equation} \rho\left(\setmat\right)
=\lim_{k\rightarrow\infty}\max_{\sigma
\in\left\{  1,...,m\right\}  ^{k}}\left\Vert A_{\sigma_{k}}...A_{\sigma_{2}%
}A_{\sigma_{1}}\right\Vert ^{1/k}.\label{eq-def.jsr}%
\end{equation}

This quantity is independent of the norm
used in (\ref{eq-def.jsr}), and is smaller than one if and only if the system is stable.  See \cite{jungers_lncis} for a recent survey on the topic.  Even though it is known to be very hard to compute, in recent years much effort has been devoted to approximating this quantity, because of its importance in applications.  One of the most successful families of techniques to approximate it makes use of convex optimization methods, like Sum-Of-Squares, or Semidefinite Programming (\cite{JohRan_PWQ,multiple_lyap_Branicky,AAA_MS_Thesis,Roozbehani2008,daafouzbernussou,LeeD06,convex_conjugate_Lyap,Pablo_Jadbabaie_JSR_journal,protasov-jungers-blondel09}). Other methods have been proposed to tackle the stability problem (e.g. variational methods (\cite{MM11}), or iterative methods (\cite{GZalgorithm})), but a great advantage of the former methods is that they offer a simple criterion that can be checked with the help of the powerful tools available for solving convex programs, and they often come with a guaranteed accuracy.
  \\As an example, the following simple set of LMIs is probably the first one that has been proposed in the literature in order to solve the problem:
\begin{equation}\label{eq-Lyap.CQ.SDP}
\begin{array}{rll}
 A_i^TPA_i&\prec&P \quad i=1,\ldots,m.\\
P&\succ&0.
\end{array}
\end{equation}
It appears that if these equations have a solution $P,$ then the function $x^TPx$ is a common quadratic Lyapunov function, meaning that this function decreases, whatever switching signal occurs.  This proves the following folklore theorem:
\begin{thm}
If a set of matrices $\setmat\mathrel{\mathop:}=\left\{ A_{1},...,A_{m}\right\} $ is such that the Equations (\ref{eq-Lyap.CQ.SDP}) have a solution $P,$ then this set is stable.
\end{thm} 

Starting with the LMIs (\ref{eq-Lyap.CQ.SDP}), many researchers have provided other methods, based on semidefinite programming, for proving the stability of a switching system.
 In all these methods, the stability criterion consists in verifying a set of \emph{Lyapunov Inequalities,} which we now describe.  The different methods amount to write a set of equations, which are parameterized by the values of the entries of the matrices in $\setmat.$ If these equations have a solution, then it implies that the set $\setmat$ is stable.

\begin{defi}
We call a \emph{Lyapunov function} any continuous, positive, and homogeneous function
$V(x):\mathbb{R}^n\rightarrow\mathbb{R}.$\end{defi}
 
\begin{defi}
Given a switching system of the shape (\ref{eq-switching}), a \emph{Lyapunov Inequality} is a quantified inequality of the shape:
\begin{equation} \forall x\in \re^n, V_i(Ax)\leq V_j(x),  \end{equation}
where the functions $V_i,V_j$ are Lyapunov functions, and $A$ is a particular product of matrices in $\setmat.$
\end{defi}
 For instance, the relations (\ref{eq-Lyap.CQ.SDP}) represent a Lyapunov inequality, because of the well known property $$P\succeq 0\quad \Leftrightarrow \quad \forall x, x^TPx \geq 0. $$  They represent the fact that the ellipsoid corresponding to the matrix $P$ is mapped into itself. We call such Lyapunov inequalities with SDP matrices and semidefinite inequalities \emph{ellipsoidal Lyapunov inequalities.}

\begin{remark}\label{remark-approx} It is important to note that the utility of such LMIs goes further than the simple stability criterion: by applying them to the \emph{scaled} set of matrices \begin{equation}\label{eq-scaled-set} \setmat/\gamma =\{ A/\gamma: A\in \setmat\}, \end{equation}  one can actually derive an upper bound $\gamma^*$ on the joint spectral radius, thanks to the homogeneity of the Definition (\ref{eq-def.jsr}): Take $\gamma^*$ the minimum $\gamma$ such that the scaled set (\ref{eq-scaled-set}) is stable.\\
This allows to provide an estimate of the quality of performance of a particular set of LMIs, as the maximal real number $r$ such that for any set of matrices $$r\gamma^*\leq \rho \leq \gamma^*. $$
In particular, it is
known (\cite{ando-shih}) that the estimate $\gamma^*$ obtained
with the set of Equations \ref{eq-Lyap.CQ.SDP} satisfies
\begin{equation}\label{eq-CQ.bound}
\frac{1}{\sqrt{n}}\gamma^* \leq\rho(\setmat)\leq\gamma^* ,
\end{equation}
where $n$ is the dimension of the matrices.  The reason for which more LMI criteria have been introduced in the literature cited above is that for some other sets of LMIs, one can prove that the value $r$ is larger than $\frac{1}{\sqrt{n}}.$
\end{remark}

In the recent paper \cite{ajprhscc11}, the authors have presented a framework in which all these methods find a common natural generalization.
The idea is that a set of Lyapunov inequalities describes a set of switching signals for which the trajectory remains stable.  Thus, a valid set of LMIs must cover all the possible switching signals, and provide a valid stability proof for all of these signals.  One contribution of \cite{ajprhscc11} is to provide a way to represent a set of LMIs with a directed labeled graph which represents all the stable switching signals (as implied by the LMIs).  Thus, one just has to check that all the possible switching signals are represented in the graph, in order to decide whether the corresponding set of LMIs is a sufficient condition for stability.  We now formally describe the construction and the result:

In the following, $\setmat$ can represent a set of matrices or the alphabet corresponding to this set of matrices. Also, for any alphabet $\setmat,$ we note $\setmat^*$ (resp. $\setmat^t$) the set of all words on this alphabet (resp. the set of words of length $t$). Finally, for a word $w\in \setmat^t,$ we note $A_w$ the product corresponding to $w:$ $A_{w_1}\dots A_{w_t}$.\\
We represent a set of Lyapunov inequalities on a directed
labeled graph $G(N, E)$. Each node of this graph corresponds to a Lyapunov function
$V_i$, and each edge is
labeled by a finite product of matrices, i.e., by a word from the
set $\setmat^*.$ 

\begin{figure}[ht]
\centering \scalebox{.3} {\includegraphics{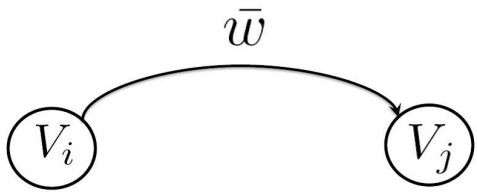}}
\caption{Graphical representation of Lyapunov inequalities. The
graph above corresponds to the Lyapunov inequality $V_j(A_wx)\leq
V_i(x)$. Here, $A_w$ can be a single matrix from $\setmat$ or
a finite product of matrices from $\setmat$.}
\label{fig-node.arc}
\end{figure}

 As illustrated in Figure~\ref{fig-node.arc}, for any word $A_w\in \setmat^*,$ and any Lyapunov inequality of the shape
\begin{equation}\label{eq-lyap.inequality.rule}
V_j(A_wx)\leq V_i(x) \quad \forall x\in\mathbb{R}^n,
\end{equation}
we add an arc going from node $i$ to node $j$ labeled with the word
$\bar w$ (the \emph{mirror} $\bar w$ of a word $w $ is the word obtained by reading $w$ starting from the end).  So, there are as many nodes in the graph as there are different functions $V_i,$ and as many arcs as there are inequalities.

\begin{figure}[ht]
\centering \scalebox{.4} {\includegraphics{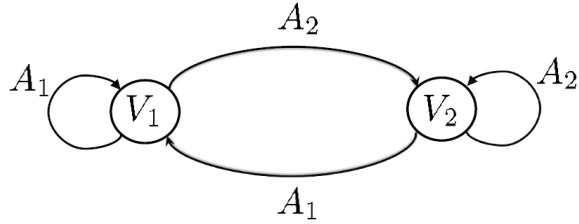}}
\caption{A graph corresponding to the LMIs in Equation (\ref{eq-hscc}).  The graph is path-complete, and as a consequence any switching system that satisfies these LMIs is stable.}
\label{fig-hscc}
\end{figure}

The reason for this construction is that there is a direct way of checking on $G$ whether the set of Lyapunov inequalities implies the stability.  Before to present it, we need a last definition:

\begin{defi}\label{def-path-complete}
Given a directed graph $G(N,E)$ whose arcs are labeled with words
from the set $\setmat^*$, we say that the graph is
\emph{path-complete}, if for any finite word $w_1\dots w_k$ of any length $k$ (i.e., for all words in
$\setmat^*$), there is a directed path in $G$ such that the word obtained by concatenating the labels of the edges on this path contains the word $w_1\dots w_k$ as a subword.
\end{defi}
We are now able to state the criterion for validity of a set of LMIs:

\begin{thm}\label{thm-path.complete.implies.stability}(\cite{ajprhscc11})
Consider a set of Lyapunov inequalities with $m$ different labels, and its corresponding graph $G(V,E).$
If $G$ is path-complete, then, the Lyapunov inequalities are a valid criterion for stability, i.e., for any finite set of matrices
$\setmat=\{A_1,\ldots,A_m\}$ which satisfies these inequalities, the corresponding switching system (\ref{eq-switching}) is stable. \end{thm}
\begin{example}
The graph represented in Figure \ref{fig-hscc} is path-complete: one can check that every word can be read on this graph.  As a consequence, the set of Equations (\ref{eq-hscc}) is a valid condition for stability.  
\end{example}

\begin{equation}\label{eq-hscc}
\begin{array}{rll}
A_1^TP_1A_1&\prec&P_1 \\
A_1^TP_1A_1&\prec&P_2 \\
A_2^TP_2A_2&\prec&P_1 \\
A_2^TP_2A_2&\prec&P_2 \\
P_1,P_2&\succ&0.
\end{array}
\end{equation}

In this paper, we investigate the converse direction of Theorem \ref{thm-path.complete.implies.stability}, and answer the question ``Are there other sets of LMIs, which do not correspond to path-complete graphs, but are sufficient conditions for stability?'' We provide a negative answer to this question. Thus, we characterize the sets of LMIs that are a sufficient condition for stability.  Of course, by Remark \ref{remark-approx} above, we not only characterize the Lyapunov inequalities that allow to prove stability, but we also characterize the valid LMIs which allow to approximate the joint spectral radius.  Another motivation for studying LMI criteria for stability is that it appeared recently that much more quantities that are relevant for the asymptotic behaviour of switching systems can also be approximated thanks to LMIs.  This is the case for instance of the Lyapunov exponent (\cite{protasov-jungers-lyap}), the p-radius (\cite{jungers-protasov-pradius}),...

Thus, we need to show that for any non-path-complete graph, there exists a set of matrices which is not stable, but yet satisfies the corresponding equations.  This is not an easy task a priori, because we need to implicitly construct a counterexample, without knowning the graph, but just with the information that it is not path-complete.  Moreover, we not only have to construct the unstable set of matrices which is a counterexample, but we need to implicitly build the solution $\{P_i\}$ to the Lyapunov inequalities, in order to show that the set satisfies the Lyapunov inequalities.

We split the proof in two steps: we first study a particular case of non-path-complete graphs: For these graphs there are only two different characters (i.e. two different matrices in the set); there is only one node; and $2^l-1$ self loops, each one with a different word of length $l$ ($l$ is arbitrary).  The proof is simpler for this particular case, and we feel it gives a fair intuition on the reasoning.\\
The rest of the paper is as follows: in Section \ref{section-preliminaries} we first present the basic construction which lies at the core of our proofs, and then we present the proof for the particular case.  In Section \ref{section-main} we prove our result in its full generality.  We then show that it implies that recognizing if a set of LMIs is a valid criterion for stability is PSPACE-complete. In Section \ref{section-conclusion} we conclude and point out some possible further work.

\section{Proof of a particular case}\label{section-preliminaries}

\subsection{The construction}

 We restrict ourselves to sets of two matrices for the sake of clarity and conciseness.  Recall that we want to prove that if a graph is not path-complete, it is not a valid criterion for stability, meaning that there must exist a set of matrices that satisfies the corresponding Lyapunov inequalities, but yet, is not stable.  Our goal in this subsection is to describe a simple construction that will allow us to build such a set of matrices in our main theorems.
If a graph is not path-complete, there is a certain word $w$ which cannot be read on the graph.  Let us fix $n=|w|+1,$ i.e., $n$ is the length of this word plus one ($n$ will be the dimension of our matrices). We propose a simple construction of a set of matrices such that all long products which are not zero must contain the product $A_w.$

\begin{defi}
Given a word $w\in \{1,2\}^*,$ we call $\setmat_w$ the set of $\{0,1\}$-matrices $\{A_1,A_2\}$ such that the $(i,j)$ entry of $A_{l}$ is equal to one if and only if
\begin{itemize}
\item $j=i+1 \mod n, \mbox{ and } w_i=l,$ for $1\leq i\leq n-1,$ 
\item { or } $(i,j)=(n,1)$ and $l=1.$
\end{itemize}\end{defi}

\begin{figure}
\centering \scalebox{0.4} {\includegraphics{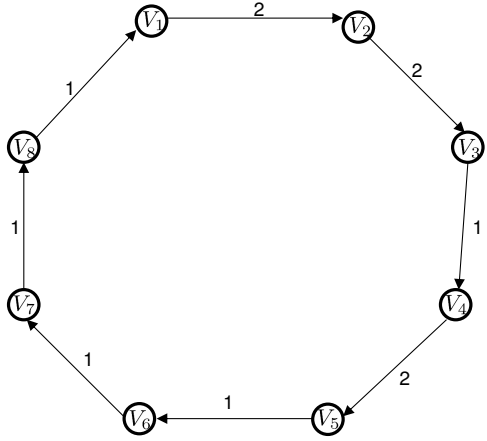}}
\caption{Graphical representation of the construction of the set of matrices $\setmat_{2212111}:$ the edges with label $1$ represent the matrix $A_1$ (i.e. $A_1$ is the adjacency matrix of the subgraph with edges labeled with a ``$1$''), and the edges with label $2$ represent the matrix $A_2.$}
\label{fig-sigmaomega}
\end{figure}
More clearly, $\setmat_w$ is the only set of binary matrices whose sum is the adjacency matrix of the cycle on $n$ nodes, and such that for all $i\in[1:n-1],$ the $i$th edge of this cycle is in the graph corresponding to $A_{w_i},$ the last edge being in the graph corresponding to $A_1.$  Figure \ref{fig-sigmaomega} provides a visual representation of the set $\setmat_w.$

\subsection{Proof of the particular case}\label{subsection-particular}
\begin{figure}
\centering \scalebox{0.3} {\includegraphics{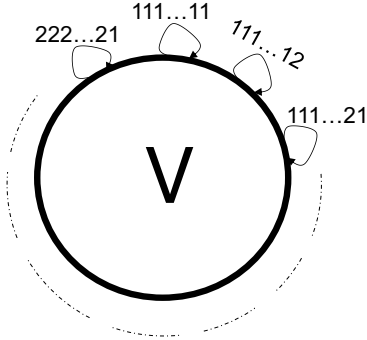}}
\caption{The graph corresponding to the particular case at stake in Subsection \ref{subsection-particular}: a single node with $2^l-1$ self loops labeled with words of length $l.$  This graph is not path-complete because the self loop $22\dots 2$ is missing.}
\label{fig-particular}
\end{figure}

We now prove a particular case of our main result: we restrict our attention to ellipsoidal Lyapunov inequalities, and to graphs with a single node with $2^{l}-1$ self-loops labeled with words of length $l.$  That is, the Lyapunov inequalities express the constraints that all but one of the products of length $l$ leave a particular ellipsoid invariant.  The corresponding graph is depicted in Figure \ref{fig-particular}, and the Lyapunov inequalities are of the shape (here we have taken $w=22\dots 2$ as the missing word):
\begin{equation}\label{eq-particular}
\begin{array}{rll}
 (A_1\dots A_1)^TP(A_1\dots A_1)&\prec&P \\
 (A_1\dots A_1A_2)^TP(A_1\dots A_1A_2)&\prec&P \\
 \dots & \dots &\dots\\
 (A_2\dots A_2A_1)^TP(A_2\dots A_2A_1)&\prec&P \\
P&\succ&0. \\
\end{array}
\end{equation}

We will make use of a well-known result from the seminal paper of \cite{ando-shih}:

\begin{thm}\label{thm-ando-shih}(\cite{ando-shih})
Let $\setmat\subset \re^{n\times n}$ be a set of matrices.  If $\rho(\setmat)<1/\sqrt{n},$ then the matrices in $\setmat$ leave invariant a common ellipsoid; that is, Equations (\ref{eq-Lyap.CQ.SDP}) have a solution. \end{thm}

We also need the easy lemma characterizing the main property of our construction $A_w.$

\begin{lem}\label{lem-subproduct}
Any nonzero product in $\setmat_w^{2n}$ contains $A_w$ as a subproduct.
\end{lem}
{\it Proof.}
Since the matrices have binary entries, a nonzero product corresponds to a path in the corresponding graph.  A path of length more than $n$ must contain a cycle, and there is only one cycle.  Finally, a path of length $2n$ must contain a cycle starting at node $1,$ which corresponds to $A_w.${\flushright $\Box$}

We are now in position to prove the particular case:
\begin{thm}
Let $w$ be a word of length $l$ on $2$ letters, and $\Sigma=\{1,2\}^{l}\setminus \{w\}$ be the set of all words of length $l$ except $w.$ The graph with one node and $2^l-1$ self-loops, whose labels are all the words in $\Sigma$ is not a sufficient condition for stability.  Indeed, the set $\setmat_w$ described above satisfies the corresponding LMIs, but is not stable.
\end{thm}
{\it Proof.}
We consider the above construction $\setmat_w.$  It is obvious that $\rho(\setmat_w)\geq 1$ (because $\rho(A_wA_1)=1;$ in fact $\rho(\setmat_w)= 1$ but this is not relevant for the discussion here).\\
Thus, we have to show that all products in the set $$\setmat'=\{A_{x_1}\dots A_{x_l}:x\in \{1,\dots,m\}^l,x\neq w\}$$ share a common invariant ellipsoid.  In order to do that, we will show that $\rho(\setmat') =0.$  This fact together with Theorem \ref{thm-ando-shih} implies that the system (\ref{eq-particular}) has a solution.

We {\bf claim} that any nonzero product $A$ of length $l=n-1$ only has nonzero entries of the shape $(i,i-1 \ (\mbox{mod}\, n)):$ $$A_{i,j}=1\,\rightarrow j=i-1\ \mbox{mod}\ n.$$
This is because any edge in the graph is of the shape $(v_i\rightarrow v_{i+1\  (\mbox{mod}\ n)}),$ so a path of length $n-1$ must be of the shape $v_i\rightarrow v_{i+n-1\  (\mbox{mod}\ n)}.$  Also, by construction of $\setmat_w,$ the only product $A$ of length $l$ such that $A_{1,n}=1$ is $A_w.$

Now, suppose by contradiction that there exists a long nonzero product of matrices in $\setmat':$ $A_{y_1}\dots A_{y_T}\neq 0,\ A_{y_i}\in \setmat'.$ Any nonzero entry in this matrix corresponds to a path of length $lT$ of the shape $v_{i_1}\rightarrow v_{i_1-1}\rightarrow\dots \rightarrow v_{i_1-T}$ (where an arrow represents a jump of length $l$ corresponding to a multiplication by a matrix in $\setmat'$).  Since we suppose that there are arbitrarily long products, it means that for some $j,$ $v_{i_j}=v_1$ and  $v_{i_{j+1}}=v_n,$ so that $A_w$ must be in $\setmat',$ a contradiction.~{\flushright{$\Box$}}

\section{The main result}\label{section-main}

\subsection{The proof}

Let us now consider a general non-path-complete graph, and prove that some sets of matrices satisfy the corresponding equations, but fail to have JSR smaller than one.  We will start by studying another family of Lyapunov functions.  These functions are only defined on the positive orthant but, as we will see, it is sufficient for nonnegative matrices.  We note these functions $V_p,$ where $p$ is a positive vector which defines them entirely:
\begin{equation}\label{eq-linear-norms}V_p(x)=\inf{\{\lambda: x /\lambda \leq p\}},\end{equation} where the inequality is entrywise.  This quantity is a valid norm for nonnegative vectors, and geometrically, its unit ball is simply the set $\{x=p-y:y\geq 0\}$. We call this family of Lyapunov functions \emph{\conicnorm{} Lyapunov functions.}  The following lemma provides an easy way to express the stability equations for this family of homogeneous functions, when dealing with nonnegative matrices: it allows to write the Lyapunov inequalities in terms of a Linear Program.

\begin{lem}\label{lem-linear-norms}
Let $p,p'\in \re^n_{++}$ be positive vectors, and $A\in \re^{n\times n}_+$.  Then, we have
\begin{equation}\label{eq-conicnorm-lp} \forall x\in \re^n_{+}, \, V_{p}(Ax)\leq V_{p'}(x)\quad \iff \quad Ap' \leq p ,\end{equation} where the vector inequalities are to be understood componentwise.
\end{lem}
{\it Proof.}
$\Rightarrow:$ Taking $x=p'$ in the left-hand side of (\ref{eq-conicnorm-lp}), and taking into account that $V_{p'}(p')=1,$ we obtain that $V_{p}(Ap')\leq 1,$ and then $$Ap'\leq p.$$
$\Leftarrow:$ First, remark that for any pair of nonnegative vectors $y,z\in\re_+^n$ and any nonnegative matrix $A\in \re^{n\times n}_+,$ $y\leq z$ implies $Ay\leq Az.$  Now, take any vector $x\in\re^n_{+},$ and denote $\gamma=V_{p'}(x).$  We have \begin{eqnarray} x/\gamma &\leq & p'\\ Ax/\gamma &\leq &Ap'\\&\leq &p. \end{eqnarray} Thus, $$V_p(Ax)=\inf\{\lambda: Ax/\lambda \leq p\}\leq \gamma. $$
{\flushright $\Box$}

We can now present our main result. 

\begin{thm}\label{theo-stab-implies-pc} A set of \conicnorm{} Lyapunov inequalities (like defined in (\ref{eq-linear-norms})) is a sufficient condition for stability if and only if the corresponding graph is path-complete.
\end{thm}
{\it Proof.}
The if part is exactly Theorem \ref{thm-path.complete.implies.stability}. We now prove the converse: for any non-path-complete graph, we constructively provide a set of matrices that satisfies the corresponding Lyapunov inequalities (with \conicnorm{} Lyapunov functions), but which is not stable.

{\bf The counterexample}
For a given graph $G$ which is not path-complete, there is a word $w$ that cannot be read as a subword of a sequence of labels on a path in
this graph.  We reiterate the construction $\setmat_{ \bar w}$ above with the particular word $\bar w.$  We show below that the set of equations corresponding to $G$ admits
a solution for $\setmat_{ \bar w}$ within the family of \conicnorm{} Lyapunov functions.

{\bf Explicit solution of the Lyapunov inequalities}
We have to construct a vector $p_i$ defining a norm for each node of the graph $G.$ In order to do this, we construct an \emph{auxiliary graph} $G'$ from the graph $G.$  The set of nodes of $G'$ are the couples ($N$ is the number of nodes in $G$ and $n$ is the dimension of the matrices in $\setmat_{ \bar w}$): $$V'=\{(i,l): 1\leq i\leq N, 1\leq l\leq n\} $$ (that is, each node represents a particular entry of a particular Lyapunov function $p_i$).
There is an edge in $E'$ from $(i,l)$ to $(j,l')$ if and only if \begin{enumerate} \item there is a matrix $A_k\in \setmat_{ \bar w}$ such that
\begin{equation}\label{eq-auxiliary-graph} (A_k)_{l',l}=1,\end{equation} \item \label{item-label}there is an edge from $i$ to $j$ in $G$ with label $A_k.$ \end{enumerate} We give the label $A_k$ to this edge in $G'.$

We {\bf claim} that if $G$ is not path-complete, $G'$ is acyclic.\\
Indeed, on the one hand, by (\ref{eq-auxiliary-graph}), a cycle $(i,l)\rightarrow \dots \rightarrow (i,l)$ in $G'$ implies the existence of a product of matrices in $\setmat_{ \bar w}$ such that $A_{l,l}=1.$  We can then build from the right to the left a nonzero product of length $2n$ by following this cycle (several times, if needed).  By Lemma \ref{lem-subproduct}, this implies that one can follow a path in $G'$ of the shape $$(i_i,l_1),\dots,(i_{n-1},l_{n-1}), $$ such that the sequence of labels is $\bar{\bar{w}}=w.$  \\
On the other hand, by item (\ref{item-label}) in our construction of $G',$ any such path in $G'$ corresponds to a path in $G$ with the same sequence of labels, a contradiction. 

Let us construct $\setmat_w$ and $G'=(V',E')$ as above.  It is well known that the nodes of an acyclic graph admit a renumbering $$s: \, V\rightarrow \mathbb{N}:\quad v\rightarrow s(v) $$ such that there can be a path from $v$ to $v'$ only if $s(v)<s(v')$ (see \cite{kahn62}).  We are now in position to define our nonnegative vectors $v_i:$  we assign the $l$th entry of $v_i$ to be equal to $s((i,l)).$

{\bf Proof that the construction is a valid solution.}
 Let us now prove that for all edge $i\rightarrow j$ of $G=(V,E)$ with label $A,$ $Av_i\leq v_j$ (where the inequality is entrywise).  This, together with Lemma \ref{lem-linear-norms}, proves that for all $x,$ $V_{v_j}(Ax)\leq V_{v_i}(x).$  
 
 Take any edge $(i,j)$ in $E$ with label $A,$ and take an arbitrary index $l'.$ Supposing $(Av_i)_{l'}\neq 0,$ we have a particular index $l$ such that $$(Av_i)_{l'}=(v_i)_{l}\leq (v_j)_{l'.}$$ This is because $A_{l',l}=1,$ together with $(i,j)\in E$ implies that there is an edge $((i,l)\rightarrow (j,l'))\in E'.$ Thus, $(v_i)_{l}< (v_j)_{l'},$ and the proof is complete.
{\flushright $\Box$}

We now provide an analogue of this result for Ellipsoidal Lyapunov functions.
\begin{thm}\label{thm-ellipsoid-contrex} A set of Lyapunov ellipsoidal equations is a sufficient condition for stability if and only if the corresponding graph is path-complete.
\end{thm}
{\it Proof.} 
The proof is to be found in the expanded version of this paper.
{\flushright $\Box$}

\begin{figure}[ht]
\centering \scalebox{.4} {\includegraphics{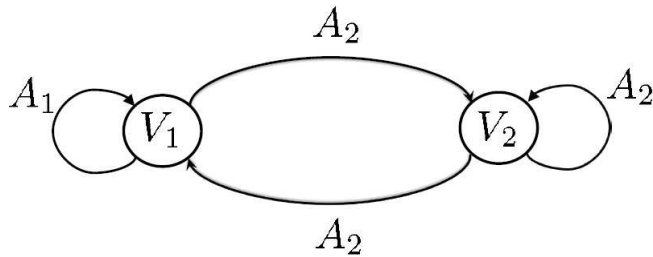}}
\caption{A graph corresponding to the LMIs in Equation (\ref{eq-hscc-wrong}).  The graph is not path-complete: one can easily check for instance that the word $A_1A_2A_1A_2\dots$ cannot be read as a path in the graph.}
\label{fig-hscc-wrong}
\end{figure}

\begin{example}
The graph represented in Figure \ref{fig-hscc-wrong} is \emph{not path-complete}: one can easily check for instance that the word $A_1A_2A_1A_2\dots$ cannot be read as a subword of a path in the graph.  As a consequence, the set of Equations (\ref{eq-hscc-wrong}) is not a valid condition for stability, even though it is very much similar to (\ref{eq-hscc}).\\
As an example, one can check that the set of matrices  $$\setmat =\left \{  \begin{pmatrix}
    -0.7   &  0.3  &   0.4\\
     0.4   &  0 &    0.8\\
    -0.7   &  0.5 &    0.7
\end{pmatrix},  
\begin{pmatrix} -0.3   & -0.95   &      0\\
    0.4   & 0.5 &   0.8\\
   -0.6       &  0   & 0.2
\end{pmatrix}\right \} $$  make (\ref{eq-hscc-wrong}) feasible, even though this set is unstable.  Indeed, $\rho(\Sigma)\geq \rho(A_1A_2A_1)^{1/3}=1.01\dots.$
\end{example}

\begin{equation}\label{eq-hscc-wrong}
\begin{array}{rll}
A_1^TP_1A_1&\prec&P_1 \\
A_2^TP_1A_2&\prec&P_2 \\
A_2^TP_2A_2&\prec&P_1 \\
A_2^TP_2A_2&\prec&P_2 \\
P_1,P_2&\succ&0.
\end{array}
\end{equation}

\subsection{PSPACE-completeness of the recognizability problem}
Our results imply that it is PSPACE-complete to recognize sets of LMIs that are valid stability criteria, as we now show.
\begin{thm}\label{thm-pspace}
Given a set of ellipsoidal Lyapunov inequalities, or a set of \conicnorm{} Lyapunov inequalities, it is PSPACE complete to decide whether they constitute a valid stability criterion.
\end{thm}
{\it Proof.} 
The proof is to be found in the expanded version of this paper.
{\flushright $\Box$}

\section{Conclusion}\label{section-conclusion}

We proved that the only sets of Lyapunov inequalities that imply stability are the ones that correspond to path-complete graphs \emph{in the case of ellipsoidal, or \conicnorm{} Lyapunov inequalities.}  
As explained above, our results are not only important for proving stability of switching systems, but also for the more general goal of approximating the joint spectral radius.

Our work leads to several interesting open questions:
It is natural to wonder whether there are other families of Lyapunov functions such that Theorem \ref{theo-stab-implies-pc} fails to be true.  That is, are there families of Lyapunov functions such that the class of valid sets of Lyapunov inequalities is larger than the path-complete ones?\\ As an example, one might look to the \emph{Complex Polytope Lyapunov functions} (see \cite{GZalgorithm,jungersprotasov09} for a study of these Lyapunov functions).  We haven't found such examples yet.  Also, the techniques analyzed here seem to be generalizable to other hybrid systems, or to the analysis of other joint spectral characteristics, like the joint spectral subradius, or the Lyapunov exponent.  Finally, our PSPACE-completeness proof does not work for graphs with two different labels, i.e., for sets of two matrices, and the complexity of path-completeness recognizability is left open in that case.

\section{Acknowledgements}
We would like to thank Marie-Pierre B\'eal and Julien Cassaigne for helpful discussions on the recognizability problem.

\def\cprime{$'$} \newcommand{\noopsort}[1]{} \newcommand{\singleletter}[1]{#1}

\newpage

\appendix
\section{Proof of Theorem \ref{thm-ellipsoid-contrex}}   \label{app-thm-ellipsoid-contrex}      

{\it Proof.} 
Again, the if part is exactly Theorem \ref{thm-path.complete.implies.stability}.

For the converse, we mimic the proof of Theorem \ref{theo-stab-implies-pc}, but we construct semidefinite matrices, and not anymore positive vectors, which satisfy the ellipsoidal Lyapunov inequalities.

Let us consider a non-path-complete graph $G(V,E),$ and the solution provided by Theorem \ref{theo-stab-implies-pc} above: we have a set of matrices $\setmat=\{A_1,A_2\},$ and nonnegative vectors $p_1,p_2,\dots,p_m$ such that $\forall e=(i,j)\in E$ with label $k,$ $$A_k v_i< v_j. $$  In fact, one can see in the proof of Theorem \ref{theo-stab-implies-pc} that the vectors $v_i$ and matrices $A_k$ have a slightly stronger property: for any matrix $A_k$ and any couple of indices $(l,l')$ such that $(A_k)_{l',l}=1,$ we have ($\vectun_l$ is the $l$th vector of the standard basis):
\begin{equation}\label{eq-vi-lifted}
 A_k \vectun_l=   \vectun_{l'},
\end{equation}
and, if moreover the edge $(i,j)$ is in $E,$
$$(A_kv_i)_{l'}=(v_i)_l<(v_j)_{l'}. $$

We are now in position to construct our solution: We take $\setmat_{\bar w}^{Tr},$ i.e., the transposes of our initial set, as the set of matrices, and we define $$P_i= \sum_l{(v_i)_l \vectun_l\vectun_l^{Tr}} ,\quad 1\leq i \leq m.$$
We claim that these semidefinite positive matrices are a solution to the Lyapunov inequalities (with ellipsoidal Lyapunov functions).  Since $\Sigma_{\bar w}^{Tr}$ is not stable, this gives us the required counterexample.\\ Indeed, for any edge $ e=(i,j)\in E$ from (\ref{eq-vi-lifted}) it is straightforward that \begin{eqnarray}A_kP_iA_k^{Tr}&=& A_k \left (\sum{(v_i)_l \vectun_l\vectun_l^{Tr}}\right ) A_k^{Tr}\\ &= &   \sum{(v_i)_l (A_k\vectun_l)(\vectun_l A_k)^{Tr}}\\&\prec&  \sum{(v_j)_{l'} \vectun_{l'}\vectun_{l'}^{Tr}}\\&&=P_j. \end{eqnarray}
{\flushright $\Box$}
\section{Proof of Theorem \ref{thm-pspace}}   \label{app-thm-pspace}  

{\it Proof.}(sketch) In the \fulllanguage, one is given a finite state automaton on a certain alphabet $\Sigma,$ and it is asked whether the language it accepts is the language $\Sigma^*$ of all the possible words. It is well known that the \fulllanguage{} is PSPACE-complete (\cite{GJ-computers-igt}).   A labeled graph corresponds in a straightforward way to a finite state automaton.  Our proof works by reduction from the \fulllanguage.   However, in order to reduce this problem to the question of recognizing whether a graph is path-complete, we must be able to transform the automaton into a new one for which all the states are starting and accepting. 

We do this by adding a new fake character $f$ in the alphabet, and connecting all accepting nodes to all starting nodes, with an edge labeled with $f.$  Now, we make all the nodes starting and accepting, and we ask whether all the words in our new alphabet $\Sigma \cup \{f\}$ are accepted in our new automaton.  That is, if the obtained graph is path-complete.  One can check that if all words can be read on the graph, then in particular, it is the case of all words starting and ending with the character $f.$  Since these words are exactly of the shape $fwf,$ where $w$ is accepted by the given automaton, our new automaton generates all the words on $\setmat \cup \{f\}$ if and only if the initial one generated all the words on $\setmat.$  \\In other words, the new graph is path-complete if and only if the initial automaton was accepting all the words on $\Sigma^*.$
{\flushright $\Box$}


\begin{thebibliography}{10}

\bibitem{AAA_MS_Thesis}
A.~A. Ahmadi.
\newblock Non-monotonic {L}yapunov functions for stability of nonlinear and
  switched systems: theory and computation, 2008.
\newblock Master's Thesis, Massachusetts Institute of Technology. Available
  from \texttt{http://dspace.mit.edu/handle/1721.1/44206}.

\bibitem{ajprhscc11}
A.~A. Ahmadi, R.~M. Jungers, P.~A. Parrilo, and M.~Roozbehani.
\newblock Analysis of the joint spectral radius via lyapunov functions on
  path-complete graphs.
\newblock In {\em Hybrid Systems: Computation and Control (HSCC'11)}, Chicago,
  2011.

\bibitem{ando-shih}
T.~Ando and M.-H. Shih.
\newblock Simultaneous contractibility.
\newblock {\em SIAM Journal on Matrix Analysis and Applications},
  19(2):487--498, 1998.

\bibitem{multiple_lyap_Branicky}
M.~S. Branicky.
\newblock Multiple {L}yapunov functions and other analysis tools for switched
  and hybrid systems.
\newblock {\em IEEE Transactions on Automatic Control}, 43(4):475--482, 1998.

\bibitem{cresp}
V.~Crespi, G.~Cybenko, and G.~Jiang.
\newblock The theory of trackability with applications to sensor networks.
\newblock {\em ACM Transactions on Sensor Networks}, 4(3):1--42, 2008.

\bibitem{daafouzbernussou}
J.~Daafouz and J.~Bernussou.
\newblock Parameter dependent {L}yapunov functions for discrete time systems
  with time varying parametric uncertainties.
\newblock {\em Systems and Control Letters}, 43(5):355--359, 2001.

\bibitem{GJ-computers-igt}
M.~R. Garey and D.~S. Johnson.
\newblock {\em Computers and Intractability; A Guide to the Theory of
  NP-Completeness}.
\newblock W. H. Freeman \& Co., New York, NY, 1990.

\bibitem{convex_conjugate_Lyap}
R.~Goebel, A.~R. Teel, T.~Hu, and Z.~Lin.
\newblock Conjugate convex {L}yapunov functions for dual linear differential
  inclusions.
\newblock {\em IEEE Transactions on Automatic Control}, 51(4):661--666, 2006.

\bibitem{GZalgorithm}
N.~Guglielmi and M.~Zennaro.
\newblock An algorithm for finding extremal polytope norms of matrix families.
\newblock {\em Linear Algebra and its Applications}, 428:2265--2282, 2008.

\bibitem{hmcb10}
E.~Hernandez-Varga, R.~Middleton, P.~Colaneri, and F.~Blanchini.
\newblock Discrete-time control for switched positive systems with application
  to mitigating viral escape.
\newblock {\em International Journal of Robust and Nonlinear Control},
  21:1093--1111, 2011.

\bibitem{JohRan_PWQ}
M.~Johansson and A.~Rantzer.
\newblock Computation of piecewise quadratic {L}yapunov functions for hybrid
  systems.
\newblock {\em IEEE Transactions on Automatic Control}, 43(4):555--559, 1998.

\bibitem{jungers_lncis}
R.~M. Jungers.
\newblock The joint spectral radius, theory and applications.
\newblock In {\em Lecture Notes in Control and Information Sciences}, volume
  385. Springer-Verlag, Berlin, 2009.

\bibitem{jungersprotasov09}
R.~M. Jungers and V.~Yu. Protasov.
\newblock Counterexamples to the {C}{P}{E} conjecture.
\newblock {\em SIAM Journal on Matrix Analysis and Applications},
  31(2):404--409, 2009.

\bibitem{jungers-protasov-pradius}
R.~M. Jungers and V.~Yu. Protasov.
\newblock Fast methods for computing the $p$-radius of matrices.
\newblock {\em SIAM Journal on Scientific Computing}, 33(3):1246--1266, 2011.

\bibitem{kahn62}
A.~B. Kahn.
\newblock Topological sorting of large networks.
\newblock {\em Communications of the ACM}, 5:558--562, 1962.

\bibitem{LeeD06}
J.~W. Lee and G.~E. Dullerud.
\newblock Uniform stabilization of discrete-time switched and {M}arkovian jump
  linear systems.
\newblock {\em Automatica}, 42(2):205--218, 2006.

\bibitem{MM11}
T.~Monovich and M.~Margaliot.
\newblock Analysis of discrete-time linear switched systems: a variational
  approach.
\newblock {\em SIAM Journal on Control and Optimization}, 49:808--829, 2011.

\bibitem{Pablo_Jadbabaie_JSR_journal}
P.~A. Parrilo and A.~Jadbabaie.
\newblock Approximation of the joint spectral radius using sum of squares.
\newblock {\em Linear Algebra and its Applications}, 428(10):2385--2402, 2008.

\bibitem{protasov-jungers-lyap}
V.~Yu. Protasov and R.~M. Jungers.
\newblock Convex optimization methods for computing the lyapunov exponent of
  matrices.
\newblock {\em Preprint}, 2011.
\newblock {A}vailable online at
  \url{http://perso.uclouvain.be/raphael.jungers/publis_dispo/protasov_jungers%
_lyapunov.pdf}.

\bibitem{protasov-jungers-blondel09}
V.~Yu. Protasov, R.~M. Jungers, and V.~D. Blondel.
\newblock Joint spectral characteristics of matrices: a conic programming
  approach.
\newblock {\em SIAM Journal on Matrix Analysis and Applications},
  31(4):2146--2162, 2010.

\bibitem{Roozbehani2008}
M.~Roozbehani, A.~Megretski, E.~Frazzoli, and E.~Feron.
\newblock Distributed {L}yapunov functions in analysis of graph models of
  software.
\newblock {\em Springer Lecture Notes in Computer Science}, 4981:443--456,
  2008.

\bibitem{pappas-multihop}
Gera Weiss, Alessandro D'Innocenzo, Rajeev Alur, Karl~Henrik Johansson, and
  George~J. Pappas.
\newblock Robust stability of multi-hop control networks.
\newblock In {\em CDC}, pages 2210--2215, 2009.

\end{thebibliography}
\end{document}